\documentclass[12pt]{amsart}

\usepackage{amssymb}  % AMS fonts
\usepackage{fullpage}
\usepackage{array}
\usepackage{enumitem}
\usepackage{url}
\usepackage{verbatim}

\theoremstyle{plain}

\theoremstyle{remark}

\numberwithin{equation}{section}

\title{Ken Kunen: Algebraist}
\author{Michael Kinyon}
\address{Department of Mathematics \\
University of Denver \\ 2360 S Gaylord St \\
Denver, Colorado 80208 USA}
\email{\url{mkinyon@math.du.edu}}

\begin{document}
\maketitle

\section{Introduction}

Ken Kunen is justifiably best known for his work in set theory and
topology. What I would guess many of his friends and students in those
areas do not know is that Ken also did important work in algebra,
especially in quasigroup and loop theory. In fact, I think it
is not an exaggeration to say that his work in loop theory, both alone
and in collaboration, revolutionized the field. In this paper, I would like
to describe some of his accomplishments to nonspecialists. My point of
view is personal, of course, and so I will give the most attention to
those of his projects in which he collaborated with me. My hope is that
the set theorists and topologists reading this will come away with
an appreciation for what Ken was able to do in an area outside of his
direct speciality. In the interest of space, I will have to leave out
discussion of some of Ken's work, such as his paper on alternative
loop rings \cite{rings}.

Ken's approach to algebra utilized automated deduction tools
and finite model builders. At the time most of what I am going to describe
took place, the automated deduction software of choice was \textsc{OTTER},
developed by William McCune \cite{Otter}. (McCune is best known to mathematicians
for his solution to the Robbins Problem in Boolean algebras \cite{Robbins}.)
The finite model builder Ken used during the period I will discuss was SEM,
developed by J. Zhang and H. Zhang \cite{SEM}. In recent
years, these have been supplanted by other tools, such as McCune's
\textsc{Prover9}, a successor to \textsc{OTTER}, and McCune's model
builder \textsc{Mace4} \cite{Prover9}.

\section{Single Axioms for Groups}

Of course, in my choices of which parts of Ken's work to discuss here,
I do not mean to suggest that his work in set theory and topology is
entirely devoid of algebraic content. But I think it is reasonable to say
that Ken's purely algebraic work begins with his papers \cite{single1,single2}
and also \cite{single3} written in collaboration with Joan Hart.

Ken's work on single axioms for groups followed up earlier work of
McCune and Wos (see the bibliographies of \cite{single1,single2,single3}).
McCune had shown that
\[
[w((x^{-1}w)^{-1}z)][(yz)^{-1}y] = x
\]
is a single axiom for group theory in the language of one binary operation $\cdot$
(written as juxtaposition) and one unary operation ${}^{-1}$.
In \cite{single1}, Ken showed that McCune's axiom is the shortest of its
particular type, consisting of $7$ variable occurrences on the left side
and $4$ variables. Ken also showed that
\[
(yy^{-1})^{-1} [(y^{-1}z)((yx)^{-1}z)^{-1}] = x
\]
is a shortest single axiom for group theory of the type with $7$ variable
occurrences on the left and $3$ variables. The follow-up papers \cite{single2,single3}
gave single axioms for groups of exponent $4$ and odd exponent, respectively.

Already in these papers we see Ken's nascent interest in nonassociative
structures. In order to verify that certain axioms are optimal, he found
it necessary to construct nonassociative models showing that other possible
candidate single axioms did not work. In retrospect, these models turned
out to be loops. Ken told me later that at the time he was writing these
papers, he did not even know what a loop was!

It is also in the single axioms papers that we find the beginnings of what
I will call the \emph{Humanization Imperative}:

\smallskip

\noindent\emph{Proofs generated by automated deduction tools should be
humanized, if it is reasonable to do so.}

\smallskip

Disclaimer: Ken never stated this imperative in precisely this form, nor even called it an imperative.
Nevertheless, we find the idea in these early algebra papers. What follows is a bit of
context for what the imperative means.

Most papers, both then and now,
in which automated deduction tools are applied to mathematical problems
do not necessarily give much human understanding of the proof
the tools generated. Generally, the problem and the strategies for using
the tools are described and the result is announced. Sometimes, the
software's proof will be included in the paper if it is not too long; 
otherwise, it
will be left to a companion website or available from the author upon request.
Of course, such proofs are usually not conceptual,
particularly those involving convoluted equational reasoning. A human
reading one of these papers is unlikely to come away with any understanding
beyond a sense that a proof exists. 

What Ken began to articulate in the single axiom papers is that it is
incumbent upon those of us who use automated deduction tools to try to
find conceptual, humanly understandable proofs, whenever doing so is
reasonable. Here is an explicit statement of this taken from \cite{single3}:

\begin{quotation}
We\ldots used \textsc{OTTER} as a reasoning assistant\ldots but found that
by examining the output from our assistant, we could provide conceptual proofs
which a human could also understand. This conceptual understanding, in turn,
led us to discover more axioms.
\end{quotation}

The Humanization Imperative plays a role in Ken's work in quasigroup and loop theory
as well.

\section{Moufang Quasigroups}

Combinatorially, a \emph{quasigroup} $(Q,\cdot)$ is a set $Q$ with a binary
operation $\cdot$ (represented by juxtaposition) such that for each $a,b\in Q$, the equations $ax = b$
and $ya = b$ have unique solutions $x,y\in Q$. In the finite
case, the multiplication tables of quasigroups are Latin squares.
A \emph{loop} is a quasigroup with a neutral element: $1x = x1 = x$
for each $x$.
A standard reference for quasigroup and loop theory is \cite{Bruck}, and
any uncited concepts in this paper can be found there.

A somewhat more useful definition is that preferred by universal algebraists.
To them, a quasigroup has not one, but three binary operations
$\cdot$, $\backslash$, $/$ satisfying
the identities
\begin{align*}
x\backslash (xy) &= y & x(x\backslash y) &= y \\
(xy)/ y &= x & (x/y)y &= x\,.
\end{align*}

A \emph{Moufang loop} is one satisfying the Moufang identities
\begin{align*}
x((yz)x) &= (xy)(zx) & (x(yz))x &= (xy)(zx) \\
x(y(xz)) &= ((xy)x)z & ((xy)z)y &= x(y(zy))\,.
\end{align*}
These identities turn out to be equivalent in loops; that is, if
a loop satisfies any one of them, then it satisfies all four.

A classical result in quasigroup theory is that an associative quasigroup
is a group. This is equivalent to proving that an associative quasigroup
has a neutral element, that is, it is a loop. Equationally, this can be
expressed in quasigroups by the identity $x\backslash x = y / y$. The
proof is easy: Start with
$xy = (x(x\backslash x))y = x((x\backslash x)y)$, then cancel $x$'s on
the left to get $y = (x\backslash x)y$, from which the desired identity
follows immediately.

In \cite{moufang}, Ken extended this result by showing that a quasigroup satisfying any one of
the Moufang identities is a loop, that is, must have a neutral element.
It follows that the Moufang identities are not only equivalent in loops,
they are, in fact, equivalent in quasigroups. In \cite{laws}, he extended this
to other classes of quasigroups and loops satisfying so-called identities of
Bol-Moufang type. For example, \emph{extra loops} are those satisfying the
identities
\[
(x(yz))y = (xy)(zy) \qquad ((xy)z)x = x(y(zx)) \qquad (yz)(yx) = y((zy)x)\,.
\]
These identities are equivalent in loops, and once again, Ken showed they
are equivalent in quasigroups.

\section{Conjugacy closed loops and G-loops}

For each $x$ in a loop $Q$, define the \emph{left} and \emph{right translation} maps
$L_x : Q\to Q$ and $R_x : Q\to Q$ by $yL_x = xy$ and $yR_x = yx$. (Here we follow
the convention that permutations act on the right of their arguments, the convention
Ken followed in his papers.) A loop $Q$ is said to be \emph{conjugacy closed}
(or a \emph{CC-loop}, for short) if the sets $\{ L_x\ |\ x\in Q\}$ and
$\{ R_x\ |\ x\in Q\}$ are each closed under conjugation. Thus for each $x,y\in Q$,
there exists $z\in Q$ such that $L_x^{-1} L_y L_x = L_z$ and similarly for the
right translations. These conditions can be easily expressed purely equationally,
thus making them ripe for exploration using automated deduction.

A (principal) \emph{isotope} of a loop $(Q,\cdot)$ is the same set $Q$ with
another loop operation $\circ$ defined by $x\circ y = (x/a)(b\backslash y)$
where $a,b$ are fixed elements of $Q$. Isotopes are of no interest to group
theorists: two groups are isotopic if and only if they are isomorphic. However,
this property does not characterize groups. There exist nonassociative loops
which are isomorphic to all of their isotopes; such loops are known as
\emph{G-loops}. (The term is imported from the Russian literature; ``G'' is probably
short for ``grouplike''.) CC-loops, for instance, are G-loops, but there are
other types as well.

CC-loops were introduced in the west by Goodaire and Robinson \cite{GR},
but went somewhat unexplored for about 17 years. Ken revived
interest in these loops in \cite{cc}. That paper gives a
detailed structure theory (partially subsumed by later work), and
also constructs examples of CC-loops, particularly ``small'' ones
(\emph{e.g.}, of order $pq$, \emph{etc.}).

In that same paper, he also showed that the equational theory of
CC-loops does not extend to all G-loops, because any identity holding
in all G-loops holds in all loops. He did this by constructing a
``universal'' G-loop which contains isomorphic copies of all finite loops.

In the follow-up paper \cite{gloops} (which actually appeared first,
because of journal backlogs), Ken studied G-loops in more detail and
showed that there are no nonassociative examples of order $3q$ where
$q$ is prime and $3$ does not divide $q-1$.

\section{Diassociative A-loops}

Now we turn to the more personal side of the story. First, a bit of
mathematical background. The \emph{multiplication group} $\mathrm{Mlt}(Q)$
of a loop $Q$ is the group generated by all left and right translation
maps. Its stabilizer of the neutral element is called the \emph{inner
mapping group} $\mathrm{Inn}(Q)$. This is a natural generalization of
the inner automorphism group of group theory. The inner mapping group
is generated by all mappings of the forms $R_x L_x^{-1}$ (measures of
noncommutativity, just as we would expect),
$R_x R_y R_{xy}^{-1}$ and $L_x L_y L_{yx}^{-1}$ (measures of
nonassociativity).

Unlike the group case, inner mappings of a loop need not act as
automorphisms of the loop. A loop for which this does occur
is called an \emph{A-loop}, a notion introduced by Bruck and Paige
in 1956 \cite{BP}. (These are now sometimes
called \emph{automorphic loops} to avoid the excessive acronymization
that has plagued the field.) Every commutative Moufang loop is
an A-loop, for instance.

A loop $Q$ is \emph{diassociative} if given any $x,y\in Q$, the
subloop generated by $x$ and $y$ is a group. Every Moufang loop
is diassociative; this is, in fact, a weak form of what is usually
known as Moufang's theorem.

Bruck and Paige studied diassociative A-loops, and noted that they
seem to have many Moufang-like properties. They did not explicitly
state a conjecture, but the subtext of the paper makes it clear that
they were conjecturing that every diassociative A-loop is a Moufang
loop. J. M. Osborn showed that this conjecture is true for
\emph{commutative} A-loops \cite{Osborn}. Osborn's paper is
a triumph of human equational reasoning.

The general conjecture sat dormant for many years, and it is
probably fair to say that most loop theorists forgot about A-loops.
Someone who was aware of it was Tom\'a\v{s} Kepka of Charles University
in Prague. He suggested it as an interesting problem to J. D. Phillips
(then of St. Mary's College in Moraga, now at Northern Michigan University).

In 1999, J. D. and Hala Pflugfelder and I were having an email conversation about every loop-theoretic topic under the sun. J. D. told me about the Bruck-Paige
problem, and the two of us discussed it at length, cc'ing Hala each time.
Finally, Hala suggested that we
should ask Ken Kunen. J. D.'s reaction was ``the set theorist?'' I had
not had a set theory class when I was a graduate student, so my reaction
was ``who?'' Hala explained to us that Ken had been getting some very
interesting results using computers.

Somewhat brashly, I sent an email message to Ken describing the problem
to him. Because of the explicit generators of $\mathrm{Inn}(Q)$ mentioned
above, the problem has an explicit equational form. Thus it can
be formulated as a problem for automated deduction tools. Ken's email reaction
was classic: 
{\tt  Interesting...} .
I later learned that ``Interesting\ldots'' (with the ellipsis)
was a clue that I had hooked him with a problem.

Ken set up the problem for \textsc{OTTER} and started working on it over the course
of a few days. (At the time, neither J. D. nor I knew how to use OTTER.)
About four or five days after I originally posed the problem to him,
J. D. and I received an email message with roughly the following header:\\
\vbox{
\begin{verbatim}
To: Michael Kinyon, J. D. Phillips
From: Ken Kunen
Subj: An ugly proof
\end{verbatim}
}
\noindent And so it was. \textsc{OTTER} had confirmed the Bruck-Paige conjecture,
and the proof was indeed quite ugly.

By this time, I had caught up with reading Ken's earlier loop theory
papers, and was familiar with what I have been calling here the
Humanization Imperative. So when J. D. asked me privately what the protocol was, I told
him that now our job was to translate the proof into human form.

Ken and I actually made
the mistake of starting at the front of the proof and working forward.
J. D. started at the end of the proof and worked backward, eventually
reaching a point where the proof relied on results which were already
in the Bruck-Paige paper. After a couple of rounds of polishing, the
paper was ready. Although most of the argument is motivated, at the
very center of the proof
is an interesting piece of equational reasoning which, to me, at least,
does not seem to be something which could have been discovered by
humans. It is nevertheless easy for a human to verify, unlike the
original \textsc{OTTER} output.

We submitted the paper \cite{diassoc} to the \emph{Proceedings
of the AMS}, the same journal that had published Osborn's paper
some 41 years earlier. It was accepted for publication five days
after we submitted it!

\section{The KK(P) Collaboration}

From that point on, Ken's work in loop theory was joint with me and
J. D. or just me. In \cite{general}, we extended Moufang's theorem
to a very wide class of loops which include both Moufang loops and
Steiner loops (which arise from Steiner triple systems) as special
cases. In this case, we did not obtain the main result directly from
\textsc{OTTER}, but we used \textsc{OTTER} to verify enough instances
of the result that we were able to construct an inductive argument.

In \cite{dicc}, we returned to CC-loops. Quite by accident, I discovered
that before Goodaire and Robinson introduced them to the western world,
they were already known in the Russian literature thanks to work of
Soikis \cite{Soikis}. In particular, a key open question in the
Goodaire-Robinson paper, restated in Ken's CC-loops paper, was solved
in 1991 by A. Basarab \cite{Basarab}: the factor loop of a CC-loop
by its nucleus is an abelian group. (Informally, the \emph{nucleus} of a loop
is the set of all elements that associate with all elements of the loop.)
It is interesting to speculate how Ken's own CC-loop paper would have
changed had he known about Basarab's earlier result. In any case,
we used Basarab's result along with Ken's structure theory to study
particular instances of diassociativity inside of CC-loops.

CC-loops which are fully diassociative are precisely the extra loops
studied earlier. Ken and I worked on these before, during and after
we met in person at the LOOPS conference in Prague in 2003. We
developed a very detailed structure theory for them \cite{extra}. Later, we
worked out a similar theory for the more general class of power-associative CC-loops
\cite{pacc}.
(A loop is \emph{power-associative} if each $1$-generated subloop is a group.)
Even now, loop theorists have told me how surprising they find the
main result of that paper to be: the factor loop of a power-associative, CC-loop
by its nucleus is an abelian group of exponent $12$.

\section{Epilogue}

After \cite{pacc}, Ken mostly went back to his set-theoretic and topological
roots. (We have a couple of other projects languishing on our back burners.)
In the meantime, both J. D. and I had started using \textsc{OTTER}
and its successor \textsc{Prover9} in our own loop theory work. And it is
actually this for which I am the most grateful to Ken. He was very patient
with all of my stupid questions about why certain input files were wrong,
what various output files meant, and so on.

I hope I have convinced the reader of the significance of Ken's work in
quasigroup and loop theory. He introduced a new methodological
paradigm: not only can and should we use automated deduction tools to assist
us, but we should also work to understand what those tools tell us.
For this, the field owes Ken a great debt.

I congratulate Ken on his change of employment status, and wish him all
the best in the years to come.

\end{document}